\documentclass{article}

\usepackage[
  lang = american, 
 paper = hardcover, openaccess,
]{ems-ecm-modified2}

\numberwithin{equation}{section}

\numberwithin{equation}{section}

 \usepackage{url} 
 \usepackage{framed}
\input xypic

\def \bfo {\begin {displaymath} }
\def \efo {\end {displaymath} }
\def \beq {\begin {eqnarray}}
\def \eeq {\end {eqnarray}}
\def \ba {\begin {eqnarray*}}
\def \ea {\end  {eqnarray*}}
\def \R {{\mathbb R}}

\def\D{{\cal D}}

\def \mbeq {\begin {eqnarray}}
\def \meeq {\end {eqnarray}}

\def\L{{\cal L}}

\def \bfo {\begin {displaymath} } 
\def \efo {\end {displaymath} } 

\def \beq {\begin {eqnarray}}
\def \eeq {\end {eqnarray}}
\def \ba {\begin {eqnarray*}}
\def \ea {\end  {eqnarray*}}

\def \F {{\cal F}}

\newcommand{\norm}[1]{\left\|#1 \right\|}

\def \Z {{\mathbb{Z}}} 
\def \C {{\mathbb{C}}} 

\def \H2s {H^{s+1}_0(\partial \M\times [0,T/2])}

\def \diam {\hbox{diam }} 
\def \dist {\hbox{dist}} 
 
\def \det {\hbox{det}\,} 
\def\bra{\langle} 
\def\cet{\rangle} 
 
\def \e {\varepsilon} 
\def \g {\tau}

\def \a {\alpha}

\def \u {\underline }

\def \pa0 {\partial _0}

\def \p {\partial}

\def \e {\epsilon}

\newcommand{\M}{{\mathcal M}}

\newcommand{\eps}{{\varepsilon}}

\newcommand{\beqla}[1] {\begin {eqnarray}\label{#1}}

\def \e {\varepsilon}

\def \beq {\begin {eqnarray}}
\def \eeq {\end {eqnarray}}
\def \ba {\begin {eqnarray*}}
\def \ea {\end  {eqnarray*}}

\def \a {\alpha}

\def \bra{\langle}
\def \cet{\rangle}

\def\bra{\langle}
\def\cet{\rangle}

\def\P{\mathbb P}
\def\W{\mathbb W}

\def \bfo {\begin {displaymath} }
\def \efo {\end {displaymath} }
\def \beq {\begin {eqnarray}}
\def \eeq {\end {eqnarray}}

\def \ba {\begin {eqnarray*}}
\def \ea {\end {eqnarray*}}

\def \D {{\cal D}}

\def \W {{\cal W}}

\def \H{{\cal H}}

\newcommand{\canberemoved}[1]{}

\newcommand{\ftext}{}
\newcommand{\vtext}{}

\newcommand{\refHOX}[1]{}

\def\itext{}
\def\gtext{} 
\newcommand{\revtext}{}

\newcommand{\removedEinstein}[1]{}

\newcommand{\commented}[1]{} 
\newcommand{\arxivpreprint}{} 

\newcommand{\motivation}{} 

\newcommand{\modified}{}

\newcommand{\extension}[1]{} 
\newcommand{\MTEXT}{}

\newcommand{\observation}[1]{}
\newcommand{\generalizations}[1]{}

\def\g{{\gamma}}
\newcommand\firstpoint{{\mathcal E}}

\newcommand{\HOX}[1]{}

\newcommand{\hiddenfootnote}[1]{}

\def\P{{\mathcal P}}

\renewcommand{\H}{{\mathcal H}}

\newcommand{\E}{{\cal E}}

\newcommand{\N}{{\mathbb N}} 

\renewcommand{\L}{{\mathcal L}}

\def \bfo {\begin {eqnarray*} }
\def \efo {\end {eqnarray*} }
\def \ba {\begin {eqnarray*} }
\def \ea {\end {eqnarray*} }
\def \beq {\begin {eqnarray}}
\def \eeq {\end {eqnarray}}

\def \dim{\hbox{dim}\,}

\def \diam {\hbox{diam }}

\def \dist {\hbox{dist}}

\def \det {\hbox{det}}

\def\bra{\langle}
\def\cet{\rangle}

\def \e {\varepsilon}
\def \p {\partial}
\def \a {\alpha}

\def\M{{\mathcal U}}
\def\F{{\mathcal F}}

\def\Z{{\mathbb Z}}

\newtheorem{definition}{Definition}[section] 
\newtheorem{theorem}[definition]{Theorem}

\newcommand{\ol}{\overline}

\newcommand{\mR}{\mathbb{R}}                    
\newcommand{\mZ}{\mathbb{Z}}                    
\renewcommand{\norm}[1]{\lVert #1 \rVert}         

\DeclareMathOperator{\Id} {Id}
\begin{document}


\title{Introduction to inverse problems for non-linear partial differential equations}
\titlemark{Inverse problems for non-linear equations}



\emsauthor{1}{
	\givenname{Matti}
	\surname{Lassas}
	\mrid{603216} 
	\orcid{0000-0003-2043-3156}}{M.~Lassas} 

\Emsaffil{1}{
	\department{Department of Mathematics and Statistics}
	\organisation{University of Helsinki}
	\rorid{ https://ror.org/040af2s02} 
	\address{P.O. Box 68 (Gustaf Hallströmin katu 2b)}
	\zip{00014}
	\city{Helsinki}
	\country{Finland}
	\affemail{Matti.Lassas@helsinki.fi}}



\begin{abstract}
We consider inverse problems for non-linear hyperbolic and elliptic equations and give an introduction
to the method based on the multiple linearization, or on the construction of artificial sources, to solve these problems. 
The method is based on self-interaction of linearized waves or other solutions in the presence of non-linearities. 
Multiple linearization has successfully been used to solve inverse problems for non-linear equation
	which are still unsolved for the corresponding linear equations. 

\end{abstract}

\maketitle

	\section{Introduction}

	In this paper, we consider inverse problems for non-linear elliptic and hyperbolic partial differential equation to demonstrate the \emph{multiple linearization} method. In this method  the nonlinearity helps one to solve inverse problems for non-linear equation
	which are still unsolved for the corresponding linear equations. Our aim is give a review on recent progress on these inverse problems
    and discuss the main ideas of the proofs in simple cases.
	
	In the case where the underlying equation is linear, a standard example is the inverse problem of Calderon \cite{calderon2006inverse}. This problem arises in Electrical Impedance Tomography (EIT), a method proposed for medical and industrial imaging, where the objective is to determine the electrical conductivity or impedance in a medium by making voltage to current measurements on its boundary. 
Another example of Calderon's inverse problem is  the determination of an unknown potential $q(x)$ in a Schr\"odinger operator $\Delta +q(x)$ from boundary measurements, first solved in \cite{SyU} in dimensions $n \geq 3$. 
	
	Let us next consider analogous inverse problems under nonlinear settings. Let $\Omega \subset \mR^n$ be a bounded domain with $C^\infty$-smooth boundary, and consider the  equation 
	\[
	\Delta u + a(x,u) = 0 \text{ in $\Omega$}.
	\]
	The Dirichlet problem for this equation is related to maintaining a temperature (or concentration or population) $f$ on the boundary. The boundary measurements for such an equation, provided that it is well-posed for some class of boundary values, may be encoded by a Dirichlet-to-Neumann map $\Lambda_a$, which maps the boundary value $f$ to the flux  (i.e., the normal derivative) $\Lambda_a(f) = \p_{\nu} u|_{\p \Omega}$ of the corresponding equilibrium state across the boundary.

Since 90's, inverse problems for nonlinear elliptic equations have been widely studied. A standard method, introduced in \cite{isakov1993uniqueness} in the parabolic case, is to show that the first linearization of the nonlinear  Dirichlet-to-Neumann  map is actually the  Dirichlet-to-Neumann  map of a linear equation, and to use the theory of inverse problems for linear equations. For the semilinear Schr\"{o}dinger equation $\Delta u +a(x,u)=0$, the problem of recovering the potential $a(x,u)$ was studied in \cite{sun1996} in dimensions $n\geq 3$, and in \cite{VictorN, imanuvilovyamamoto_semilinear} in dimension $n=2$. In addition, inverse problems have been studied for quasilinear elliptic equations \cite{sun1997inverse, kang2002identification, munozuhlmann}, the degenerate elliptic $p$-Laplace equation \cite{salo2012inverse, branderetal_monotonicity_plaplace}, and the fractional semilinear Schr\"odinger equation \cite{lai2019global}. Certain Calder\'on type inverse problems for quasilinear equations on Riemannian manifolds are considered in~\cite{lassas2018poisson}, see also referenceis in \cite{Agnelli,uhlmann2009calderon}.

%
%
%
%
%
%
	
	Let us introduce  the mathematical setting to study canonical examples of non-linear elliptic equations. We will denote by $(M,g)$ a compact Riemannian manifold with $C^\infty$ boundary $\p M$, where $\dim(M)=n$, $n\geq 2$.  Let $q \in C^{\infty}(M)$. We will consider semilinear elliptic equations of the form 
	\beq\label{Main semilinear equation}
	\begin{cases}
	\Delta_g u +q u^m =0 & \text{ in }M, \\
	u=f & \text{ on }\p M,
	\end{cases}
	\eeq
	where 
	$
	m\in \N \text{ and } m\geq 2.
	$
    Here $\Delta_g$ is the Laplace-Beltrami operator, given in local coordinates by 
	\[
	\Delta_g u=\frac{1}{\det(g)^{1/2}}\sum_{j,k=1}^n\frac{\p}{\p x^j}\left( \det(g)^{1/2} g^{jk}\frac{\p u}{\p x^k}\right),
	\]
where $g=(g_{jk}(x))$ and $g^{-1}=(g^{jk}(x))$. As usual, we denote the Euclidean Laplace operator,
corresponding to the Euclidean metric $g^{(e)}=\delta_{jk}$ in $\R^n$ by $\Delta.$

	The Dirichlet problem for \eqref{Main semilinear equation} with a boundary condition
	$u|_{\p \Omega}=f$ has a unique small solution $u$ for sufficiently small boundary data $f\in C^s(\p M)$, where $s>2$ with $s\notin \N$. More precisely this means that there is $\delta>0$ such that whenever $\norm{f}_{C^s(\p M)}\leq \delta$, there is a unique solution $u$ to \eqref{Main semilinear equation} satisfying $\|u\|_{C^s(M)}\le c\delta$. We will call $u$ the unique small solution and denote it by $u_f$ to indicate the boundary value. Here $C^s$ is the standard H\"older space for $s>2$ with $s\notin \N$ (often written as $C^{k,\alpha}$ if $s = k+\alpha$ where $k \in \mZ$ and $0 < \alpha < 1$). Hence, the  Dirichlet-to-Neumann  map is defined by using the unique small solution in a following way:
	\ba
	 & &\Lambda_{M,g,q}: \{f\in C^s(\p M)\mid \norm{f}_{C^s(\p M)}\leq \delta\}\to C^{s-1}(\p M), \\
	 & & \Lambda_{M,g,q}:f \mapsto \p_{\nu} u_f|_{\p M},
	\ea
	where $\p_{\nu}$ denotes the normal derivative on the boundary $\p M$. In what follows, we denote the  Dirichlet-to-Neumann  map 
    by  $\Lambda_{M,g}$ to 
    when $q=0$. When $M=\Omega \subset \R^n$ and $g$ is the identity matrix, we denote the  Dirichlet-to-Neumann  map by $\Lambda_{q}$.
	
	
	Next, we consider the case when $\Omega \subset \R^n$ for $n\geq 2$ and $\Delta_g$ is the Euclidean Laplacian.
	
	\begin{theorem}[Reconstruction of non-linear term \cite{lassas2019nonlinear}]\label{main thm q}
		Let $n\geq 2$, and let $\Omega\subset \R^n$ be a bounded domain with $C^\infty$ boundary $\p \Omega$. Let $q_1,q_2\in C^\infty(\ol{\Omega})$. Assume the  Dirichlet-to-Neumann  maps $\Lambda_{q_j}$ for the equation
        \eqref{Main semilinear equation} with potentials $q_j$,
$j=1,2$ satisfy 
		\[
		\Lambda_{q_1}(f)=\Lambda_{q_2}(f) 
		\]
		for all $f\in C^s(\p \Omega)$ with $\norm{f}_{C^s(\p M)}<\delta$, where $\delta>0$ is any sufficiently small number. Then $q_1=q_2$ in $\Omega$.
	\end{theorem}

    In the pioneering studies of the problem, the first and the second order linearizations of the nonlinear  Dirichlet-to-Neumann  map has been used in the works \cite{sun1996, sun1997inverse} related to nonlinear equations with matrix coefficients and in \cite{kang2002identification} for a nonlinear conductivity equation (see also \cite{carsteanakamuravashisth}).

    Next, we present the outline of the proof of Theorem \ref{main thm q}.
	For the equation \eqref{Main semilinear equation} with quadratic nonlinearity, that is, $m=2$, the first linearization of the nonlinear  Dirichlet-to-Neumann  map $\Lambda_q$, linearized at the zero boundary value, is just the  Dirichlet-to-Neumann  map for the standard Laplace equation:
	\[
	D|_{f=0}\Lambda_q(f): C^s(\p \Omega) \to C^{s-1}(\p \Omega), \ \ f \mapsto \p_{\nu} v_f|_{\p \Omega},
	\]%
	where $s>1$ is not an integer, $v_f$ is the unique solution of $\Delta v_f = 0$ in $\Omega$ with $v_f|_{\p \Omega} = f$. Thus the first linearization does not carry any information about the unknown potential $q$. For a quadratic nonlinearity the \emph{second linearization} of the map $f\to \Lambda_q(f)$, denoted by
    $D^2|_{f=0} \Lambda_q(f)[\cdot,\cdot]$, is a symmetric bilinear form from $C^s(\p \Omega) \times C^s(\p \Omega)$ to $C^{s-1}(\p \Omega)$. 
	The second order linearization is given by
	\beq \label{formula A}
			D^2|_{f=0} \Lambda_q(f)[f_1, f_2] = \p_{\eps_1} \p_{\eps_2} u_{\eps_1 f_1 + \eps_2 f_2}|_{\eps_1=\eps_2=0}\quad \text{ on } \p \Omega.
	\eeq
	That is, one considers boundary data 
	$
	f=\eps_1 f_1 + \eps_2 f_2\in C^s(\p \Omega),
	$
	where $\eps_1,\eps_2$ are sufficiently small parameters, and takes the mixed derivative 
\[
w(x)=    
\left.\frac{\p}{\p \eps_1}\frac{\p}{\p \eps_2}\right|_{\eps_1=\eps_2=0}u_{\eps_1 f_1 + \eps_2 f_2}(x)
\]
of the solutions of $u_{\eps_1 f_1 + \eps_2 f_2}(x)$ of the equation~\eqref{Main semilinear equation} with the coefficient $q$, $m=2$, 
and the boundary value $f=\eps_1 f_1 + \eps_2 f_2$.
		This  yields the equation
    \begin{equation}\label{second_deriv}
     \Delta w = -2q v_{f_1} v_{f_2},
    \end{equation}
     where 
    $v_{f_k}$ are harmonic functions, i.e.\ the solutions to the linearized equation $\Delta v = 0$
    with the Dirichlet boundary value $v_{f_2}|_{\p \Omega}=f_k$. Moreover,
    \beq \label{formula AB}
			D^2|_{f=0} \Lambda_q(f)[f_1, f_2] = \p_\nu w|_{\p \Omega}.
	\eeq
    Assume next that there are two potential functions $q_1$ and $q_2$ such that
    $$
    \Lambda_{q_1}(f)=\Lambda_{q_2}(f),\quad\hbox{for all }f\in C^s(\p M),\ \ \norm{f}_{C^s(\p M)}\leq \delta.
    $$    
    Taking the mixed  2nd order derivative of the  Dirichlet-to-Neumann  maps yields 
   $  \p_\nu w_1=\p_\nu w_2 $ on $ \p \Omega$,
  where $     \Delta w_j = -2q_j v_{f_1} v_{f_2},$ and $w_j|_{\p \Omega}=0.$
    Subtracting the equations~\eqref{second_deriv} for $q_1$ and $q_2$, and integrating the resulting equation against the harmonic function $v_{f_3}$ and using integration by parts yields the  formula 
    \beq\label{integration by parts}
    \int_{\Omega} (q_1-q_2) v_{f_1} v_{f_2} v_{f_3} \,dx 
      &=&-\frac 12  \int_{\Omega} (\Delta w_1-\Delta w_2)  v_{f_3} \,dx 
    \\ \nonumber
      &=&-\frac 12  \int_{\p\Omega} (\p_\nu w_1-\p_\nu w_2) v_{f_3} \,dS(x) 
    =0.
    \eeq
To apply this formula, let  us denote $Q=q_1-q_2$. When $k\in \R^3$, there are 
$\xi_1,\xi_2\in \C^3$ satisfying
\ba
\xi_1\,\cdotp \xi_1=0,\quad \xi_2\,\cdotp \xi_2=0,\quad \xi_1+\xi_2=ik,
\ea
where $\eta\,\cdotp \zeta=\eta_1\zeta_1+\eta_2\zeta_2+\eta_3\zeta_3$
for $\eta=(\eta_1,\eta_2,\eta_3),\zeta=(\zeta_1,\zeta_2,\zeta_3)\in \C^3.$
We define
\ba
& &v_1(x)=e^{\xi_1\cdotp\, x},\quad v_2(x)=e^{\xi_2\cdotp\, x},\quad v_3(x)=1.
\ea
 These functions satisfy $\Delta v_j=0$ and we can use formula \eqref{integration by parts}
 with ${f_j}=v_j|_{\p \Omega}$, that is, $v_{f_j}=f_j$.
The  Fourier transform  $\widehat Q(k)$  of $Q(x)$ (that is extended as a zero function in $\R^n\setminus \Omega)$ 
satisfies 
 \ba
& &\hspace{-2mm} \widehat Q(k)=
\int_\Omega Q(x)e^{ik\cdotp\, x}\, dx=
\int_\Omega Q(x)e^{\xi_1\cdotp\, x}\cdot e^{\xi_2\cdotp\, x}\cdot 1\, dx=\int_\Omega Qv_1v_2v_3\, dx=0.
\ea
Hence $Q=0$, that is, $q_1=q_2$. These arguments (with some details) prove Thm.\ \ref{main thm q}.

%
%
%
    Analogous inverse problems  for various non-linear elliptic equations using multiple linearization has been studied in \cite{CarsteaKrupchyk_nonlinear,FeizmohammadiKian,KianKrupchyk_nonlinear,Kru1,Kru2,Kru3}.

We move on to describe our next result. By using higher order linearizations we can solve the following \emph{simultaneous recovery} on a two-dimensional Riemannian surface and coefficients of equations.

	\begin{theorem}[Simultaneous recovery  of metric and potential \cite{lassas2019nonlinear}]\label{main thm 2D}
	Let  $(M_1,g_1)$ and $(M_2,g_2)$ be two compact connected manifolds with mutual $C^\infty$ smooth boundaries $\p M_1=\p M_2=:\p M$
	and $\dim (M_1)=\dim (M_2) =2$. Let $m\geq 2$, and let $\Lambda_{M_j,g_j,q_j}$ be the  Dirichlet-to-Neumann  maps   $\Lambda_{M_j, q_j}: f\mapsto \p_{\nu} u|_{\p M}$ of the problem
	\beq \label{2Deq}
	\Delta_{g_j} u +q_j u^m=0  \text{ in }M_j
	\eeq
	for $j=1,2$.
	Let $s>2$ with $s \notin \N$ and assume that 
	\[
	\Lambda_{M_1,g_1,q_1} f =\Lambda_{M_2,g_2,q_2} f \text{ on }\p M, 
	\]
	for any $f\in C^s(\p M)$ with $\norm{f}_{C^s(\p M)}\leq \delta$, where $\delta >0$ is sufficiently small. Then there exists a  diffeomorphism $\Psi:M_1\to M_2$ and a positive smooth function $\sigma$ on $M_1$ such that for $x\in M_1$ we have
		\begin{equation*}\label{J_and_sigma}
		(\sigma \Psi^*g_2)(x) = g_1(x),\quad \sigma q_1 = q_2 \circ \Psi \text{ in }M_1, 
		\end{equation*}
		with $\Psi|_{\p M}=\Id \text{ and }\sigma|_{\p M} =1$.
	\end{theorem}
	
        For the linear equation $\Delta_g u + qu = 0$, an analogous result has been proved when $M$ is a domain in $\mR^2$ with a Riemannian metric \cite{iuy_general}, when $M$ is a manifold and the potential $q$ is zero \cite{lassas2001determining}, and when the conformal class of the manifold $(M,g)$ is a priori known \cite{guillarmou2011calderon}. The recovery of properties of both the manifold and potential upto the gauge transformation \eqref{J_and_sigma} was a long-standing open question that was
        was solved in 2024 in by Carstea, Liimatainen, and Tzou, see \cite{Liimatainenetal} using much more advanced techniques that the problem for the non-linear equation needs.


Next, we consider the \emph{inverse obstacle problem} for elliptic equations with power type nonlinearities. Let $\Omega \subset \R^n$ be a bounded, connected open set with a $C^\infty$-smooth boundary $\p \Omega$. Let $D$ be an unknown cavity in $\Omega$. More precisely,  $D\subset \Omega$ is an open set with a smooth boundary $\p D$
such that $\overline{D}\subset \Omega$. Moreover, let $m\in \N$ with $m\geq 2$.
 Consider the following semilinear elliptic equation,  with an unknown coefficient $q=q(x)$ and 
 an unknown relatively compact domain $D\subset \Omega$,
\begin{eqnarray}\label{main equation_cavity}
	\begin{cases}
	\Delta u + qu^ m =0 & \text{ in }\Omega \setminus \overline{D},\\
	u =0 &  \text{ on }\p D, \\
	u =f & \text{ on }\p \Omega.
	\end{cases}
\end{eqnarray}
For $s>1$ and $s \in \N$, let $f \in C^s(\p M)$ with $\norm{f}_{C^s(\p \Omega)}< \delta$, where $\delta>0$ is any sufficiently small number. Then by the well-posedness of \eqref{main equation_cavity} \cite[Proposition 2.1]{lassas2019nonlinear} again, one can define the corresponding  Dirichlet-to-Neumann  map $\Lambda_{D,q}$ on $\p \Omega$ with
 $$
 \Lambda_{D,q} :  \{f\in C^s(\p M)\mid \norm{f}_{C^s(\p M)}\leq \delta\} \to C^{s-1}(\p \Omega),  \quad \Lambda_{D,q}: f \mapsto \p_\nu u_f |_{\p \Omega},
 $$ 
 where $s>1$ is not an integer and $u$ is the unique small solution of \eqref{main equation_cavity}.
 The inverse obstacle problem is to determine both the unknown cavity $D$ and potential $q$ from the  Dirichlet-to-Neumann  map $\Lambda_{D,q}$.


\begin{theorem}[Simultaneous recovery: Cavity and potential \cite{LLYS-obstacles}]\label{thm: Nonlinear Schiffer's problem}
	Let $\Omega \subset \R^n$ be a bounded domain for $n\geq 2$ with a $C^\infty$-smooth boundary $\p \Omega$. Let $D_1, D_2\Subset \Omega$ be nonempty open cavities with $C^\infty$-smooth boundaries and $q_j\in C^\infty(\Omega\setminus \overline{D_j})$. Let $\Lambda_{D_j, q_j}: f\mapsto \p_{\nu} u|_{\p \Omega}$, $j=1,2$ , be the  Dirichlet-to-Neumann  maps for
    the equations \eqref{main equation_cavity}
    with $D=D_j$ and $q=q_j$,
	$j=1,2$. Assume that $\Lambda_{D_1, q_1}f= \Lambda_{D_2,q_2}f$ for all  sufficiently small $f\in C^\infty(\p \Omega)$ on $\p \Omega$. Then, 
	$
	D_1 = D_2 $ and $q_1 =q_2$ in $\Omega\setminus \overline{D_1}.
	$
\end{theorem}

The simultaneous recovering of an obstacle and an unknown surrounding potential in dimensions $n\ge 3$ 
is  a long-standing open problem. This problem is closely related to the \emph{partial data} Calder\'on problem studied in \cite{kenig2007calderon,imanuvilov2010calderon}, where the boundary is a prior known but the boundary measurements are done only on a subset of the boundary. 

\section{Inverse problems for waves}

Let us consider a  linear wave equation 
\begin{equation*}
\partial_t^2 u(t,x) - c(x)^2 \Delta u(t,x) = 0\quad\hbox{in }\R\times\Omega
\end{equation*}
where $\Omega\subset \R^n$ and  $c(x)$ is the wave speed. This equation models e.g.\ acoustic waves in isotropic medium.
In inverse problems  one has access to measurements of waves (the solutions $u(t,x)$) on the boundary, or in a subset of the domain $\Omega$, and one aims to determine unknown coefficients {\ftext (e.g., $c(x)$)} in the interior of the domain.  

 We will consider \emph{anisotropic} materials, where the wave speed depends on the direction of propagation. This means that the scalar wave speed $c(x)$, where $x=(x^1,x^2,\dots,x^n)\in \Omega$, is replaced by a positive definite symmetric matrix $(g^{jk}(x))_{j,k=1}^n$. Anisotropic materials appear frequently in applications such as in  seismic imaging, where one wishes to determine the interior structure of 
{\ftext the} Earth by making various measurements of waves {\ftext on its} surface.

To model waves in anisotropic medium mathematically,
let $(N,g)$ be an  $n$-dimensional Riemannian manifold with boundary and consider the wave equation
\beq
\label{eq:problem}
& &\p_t^2 u(t,x) - \Delta_g u(t,x) = 0 \quad \text{in }(0,\infty) \times N,\\
& &\p_\nu u|_{\R_+\times \p N} = f,\quad
u|_{t = 0} =0,\quad  \p_t u|_{t=0} = 0,\nonumber
\eeq
where $\Delta_g$ is the Laplace--Beltrami operator corresponding 
to a smooth time-independent Riemannian metric $g$ on $N$. 
 We denote the unit speed geodesic of $(N,g)$ emanating from a point $(y,\eta)\in SN$  by $\alpha_{y,\eta}(t)=\exp_p(t\eta)$. 

{\ftext
The solution of (\ref{eq:problem}), corresponding to the boundary value $f$ (which is interpreted as a boundary source), is denoted by 
$u^f=u^f(t,x)$.}

Let us assume that the boundary $\p N$ is known.
The inverse problem is to reconstruct the manifold $N$ and the metric $g$ when we are given
the boundary $\p N$, as a manifold, and 
the {\it response operator}
\beq \label{1a}
  \Lambda_{N,g}:f\mapsto u^f|_{\R_+\times \p N},
\eeq 
which is also called the {\it Neumann-to-Dirichlet map}. Physically, 
$\Lambda_{N,g}f$ describes the measurement of the medium
response to any applied boundary source $f$. 

%

It {\ftext is} convenient to interpret the anisotropic wave speed $(g^{jk}(x))$ also in a domain of the Euclidean space as the inverse of a Riemannian metric, thus modelling the medium as a \emph{Riemannian manifold}.
This is due to fact that if $\Psi:N\to N$  is a diffeomorphism such that $\Psi|_{\p N}=Id$  then all boundary measurements for the metric $g$  and the pull-back metric $\Psi^*g$
coincide. Thus to prove uniqueness results for inverse problems, one has to consider properties that are invariant in diffeomorphisms and try to reconstruct those uniquely, for example,  to show that an underlying manifold structure {\ftext  can} be uniquely determined. In practice, the inverse problem in a subset of the Euclidean space is solved in two steps. The {\ftext first} is to reconstruct the underlying manifold structure. The second step is to find an embedding of the constructed manifold {\ftext in} the Euclidean space using additional a priori information, see \cite{KKL}. Below, we will concentrate on the first step.

{\ftext In 1990s,} the combination of the  {boundary control} method developed by Belishev \cite{Bel} and generalized for Riemannian manifolds in  \cite{BelKur}  and Tataru's unique continuation theorem \cite{Tataru1}  gave a solution to the inverse problem of determining the isometry type of a Riemannian manifold 
 $(N,g)$ with given 
 boundary $\p N$ and the Dirichlet-to-Neumann map $\Lambda_{N,g}$.

\HOX{Check citation style and ;}

\begin{theorem}[\cite{BelKur,Tataru1}]
\label{th:BK}
Let $(N_1,g_1)$ and $(N_2,g_2)$ be compact  Riemannian manifolds with boundary.
Assume that there is a diffeomorphism $\Phi:\p N_1\to \p N_2$
such that
\begin{equation}
\label{eq: phi pulls back Lambda bnd}
\Phi^\ast (\Lambda_{N_1,g_1} f)= \Lambda_{N_2,g_2} (\Phi^\ast f), \quad \hbox{for all } f\in C^{\infty}_0(\R_+\times \p N_1).
\end{equation}
Then $(N_1,g_1)$ and $(N_2,g_2)$ are isometric Riemannian manifolds.
\end{theorem}

Above, $\Phi^*f(t,x)=f(t,\Phi(x))$ is the pull-back of $f$ in $\Phi$, or more precisely, that of $Id_\R \times\Phi:
\R \times \p N_1\to \R \times\p N_2$.

The original method to prove Theorem \ref{th:BK} was based on combining the spectral theory for the Laplace operator and control theory for the wave equation.
{An alternative way to prove Theorem \ref{th:BK} is to show that the hyperbolic Dirichlet-to-Neumann
map $\Lambda_{N,g}$ can be used to determine a sequence of Dirichlet boundary values $f_j$
for which the corresponding waves $u^{f_j}(x,T)$ focus at  time $T>2\,\diam(N)$, see \cite{bingham2008iterative,Kirpichnikova-Dahl,Kirpichnikova-Korpela}:

\begin{theorem}[Focusing of waves \cite{Kirpichnikova-Korpela}]\label{thm: focusing} Let $T>2\,\hbox{diam}(N)$, $z\in \p N$, $s>0$,
and $y\in N$ be a point on the normal geodesic emanating from $z$, $y=\alpha_{z,\nu}(s)$ such that $\dist(y,\p N)=s$.
Using the Neumann-to-Dirichlet map  $\Lambda_{N,g}$, one can 
construct a sequence of boundary sources
$ f_j\in C^\infty_0((0,T)\times \p N)$, $j=1,2,\dots$ such that for some $C_y \not =0$ at the time $T$
\ba\label{focusing}
\lim_{j\to \infty}
\left(\begin{array}{c}
u^{f_{j}}(T,x)  \\
\p_t u^{ f_{j}}(T,x)
\end{array}\right)
=
\left(\begin{array}{c}
0\\
C_y \delta_{y}(x)  \\
\end{array}\right),\quad \hbox{for }x\in N
\ea
in sense of distributions in $\mathcal D'(N)\times \mathcal D'(N)$.
\end{theorem}

Let $G(x,t;y)$ be 
Green's function of the wave operator, that is, the wave produced by a point source at the point
$y$ at the time zero,
\ba
& &(\p_{t}^2-\Delta_g)G(t,x;y)=\delta_{(0,y)}(t,x)\quad\hbox{on }(t,x)\in   \R_+ \times N,\\
& & G|_{(t,x)\in  \R_+ \times\p  N }=0
,\quad G|_{t=0}=0,\quad \p_tG|_{t=0}=0,
\ea
where $\Delta_g$ operates in the $x$ coordinate.
Theorem \ref{thm: focusing} implies that the waves $u^{f_{j}}(t,x)$ converge in sense of distributions to $G(t-T,x;y)$ in the domain $x\in N$, $t>T$. This means that using the Dirichlet-to-Neumann map, it is possible to construct functions 
$(\Lambda_{N,g}f_j)(x,t-T)$ that converge, as $j\to \infty$, to the boundary
values  $\p_\nu G(\cdot,\cdot;y)|_{\R_+\times \p N}$ of the waves $G(\cdot,\cdot;y)$ produced by points sources
$C_y \delta_{y}(x)$ at the time $t=0$. 
The first time $t>0$ when the wave $G(\cdot,\cdot;y)$ is observed to be non-vanishing at the boundary point $x\in \p N$ is defined to be
\beq
{\bf t}_{y}(x)&=&\sup \{t\in \R;\hbox{ the interval $\{z\}\times (0,t)$ has a neighborhood}\\ \nonumber
& &\quad \quad \quad \quad \quad \hbox{ $U\subset \R\times \p N$ such that $ \p_\nu G(\cdot,\cdot,y)\big|_U =0$}\}.
\eeq
Using the finite propagation speed of waves and Tataru's unique continuation theorem \cite{Tataru1}, it turns out  that
$
{\bf t}_{y}(x)=d_g(x,y),
$
see \cite{LSaksala}.
Thus  Theorem \ref{thm: focusing} implies that  the Dirichlet-to-Neumann
map $\Lambda_{N,g}$ for the wave equation determines the functions
\beq\label{Riemannian observation times}
{\bf t}_{\alpha_{z,\nu}(s)}(x)=d_g(x,\alpha_{z,\nu}(s)),\quad \hbox{for 
$x,z\in \p N$ and $s>0$ such that dist$(x,z)=s$}.
\eeq
Later, we will return to these data, but before that we consider consequences of 
 Theorem \ref{th:BK}.
 
}

Theorem \ref{th:BK} can also be used to prove the uniqueness of several other inverse problems. By  \cite{KKLM}, the equivalence of spectral inverse problems with several different measurements, that in particular implies the following result.

\begin{theorem}[\cite{KKLM}]
\label{th:eqivalence} Let  $(N,g)$ a be compact smooth  Riemannian manifold with boundary $\p N$.
Then the boundary $\p N$ and Neumann-to-Dirichet map $\Lambda: \partial_\nu u|_{\R_+\times \p N} \mapsto u|_{\R_+\times \p N}$, for heat equation $(\p_t-\Delta_g)u=0$, or for
the Schr\"dinger  equation $(i\p_t-\Delta_g)u=0$, with vanishing initial data $u|_{t=0}=0$ determine the Neumann-to-Dirichet map for the wave equation,
and therefore, the manifold $(N,g)$ up to an isometry.
\end{theorem}

The stability of the solutions of the above inverse problems have been analyzed in \cite{AKKLT,Bao,Bosi-Stability,Yamaguchi1,Yamaguchi2,StU}.

Without making strong assumptions {\ftext about} the geometry of the manifold, the existing uniqueness results for the inverse problems for the linear hyperbolic equations
with vanishing initial data  are limited to equations {\ftext whose} coefficients  are 
{\ftext time independent or real-analytic} in time (see e.g.\ \cite{AlexakisOF,AKKLT,BelKur,Eskin,helin2016correlation,
KKL,Krupchyk-spectral,Krupchyk-energy,KOP,Lassas-O_duke}). The reason for this is that
these results are based on Tataru's unique
continuation theorem \cite{Tataru1}.
This sharp unique continuation result does not work for general wave equations {\ftext whose} coefficients are not
{\ftext real-analytic} in time, as shown by
Alinhac \cite{Alinhac}.
Solving the inverse problem of finding the Riemannian metric (up to an isometry), i.e. proving Theorem \ref{th:BK},
can be reduced to determining the metric from passive observations from
sources inside the domain, and due to this we next consider passive observations
from the point sources.

\subsection{Inverse problems with passive measurements}

\subsubsection{Lorentzian light observation sets and waves on a Riemannian manifold}

Let us consider a body in which there   {spontaneously appear} point sources that create propagating waves. 
In this subsection we assume that the velocity of the waves is independent of time and depends only the the spatial
variables and model the wave speed using a Riemannian metric.
 We consider the inverse problem where we detect such waves either outside or at the boundary of the body and aim to determine the unknown wave speed inside the body. 
 {As an example of such situation, one can consider the micro-earthquakes that appear very frequently 
near  active faults. The related inverse problem is whether the surface observations
of elastic waves produced by the micro-earthquakes can be used in the geophysical imaging
of Earth's subsurface, that is, to  determine the speed of the elastic waves in the studied volume.} 
Below, we consider a highly idealized version of the above inverse problem.
We consider the problem on an $n$  dimensional manifold $N$ with a Riemannian metric $h=h_{ij}(x)$
that corresponds to the travel time of a wave between two points (below, we will consider Lorentzian metric tensors for which we reserve the symbol $g$, and use now letter $h$ for the Riemannian metric). The Riemannian distance
of points $x,y\in N$ is denoted by $d_h(x,y)$. For simplicity
we assume that the manifold $N$  is compact and has no boundary.
Instead of considering measurements on boundary, we assume that the manifold
contains a known part $F$ and the metric is unknown outside this set, that is, in $N\setminus F$.
When  {a spontaneous point source} produces a wave at some unknown point $x\in N\setminus F$
at some unknown time $t\in \R$, the produced wave is observed at the point
$z\in F$  at time $T_{x,t}(z)=d_h(z,x)+t$. Note that these functions
are analogous to functions \eqref{Riemannian observation times} (that were defined
for the manifolds with boundary) with the source
points $x=\alpha_{z,\nu}(s)$. These observation times
at two points $z_1,z_2\in F$ determine
the {\it distance difference function}
\beq\label{def: dist diff}
D_x(z_1,z_2)=T_{x,t}(z_1)-T_{x,t}(z_2)=d_h(z_1,x)-d_h(z_2,x).
\eeq
Physically, this function corresponds to the difference of times
 when the waves produced by the point source at $(x,t)$
are observed at $z_1$ and $z_2$.
An assumption there are a large number point sources and that we do measurements over
a long time interval can be modeled by the assumption that we are given  {the set $F$ and}
the family of functions
$
\{D_x \mid x\in X\}\subset  C(F\times F),
$
where $X\subset N$ is either the whole set $N\setminus F$ or its dense subset.

Next, let $(N_1,h_1)$ and $(N_2,h_2)$ be compact and connected Riemannian manifolds without boundary. 
Let $d_{h_j}(x,y)$ denote the Riemannian distance of points $x,y\in N_j$, $j=1,2$. 
For $x\in N_j$ and $z_1,z_2\in F_j$ we denote $D^j_x(z_1,z_2)=d_{h_j}(x,z_1)-d_{h_j}(x,z_2)$.


%

 To prove
the uniqueness of the passive imaging inverse problem, we assume the following:
%
%
%
%
%
\begin{enumerate}
\item  [(i)] The map $\phi:F_1 \rightarrow F_2,$ is an isometry, that is, $d_{h_1}(z,w)=d_{h_2}(\phi(z),\phi(w))$ for all $z,w\in F_1$.
\item  [(ii)] The collections $\D_j(N_j)=
\{D^j_x(\cdot,\cdot);\ x\in N_j\}\subset C(F_j\times F_j)$ with $j=1,2$ are equivalent, that is,
$\{D^1_x(\cdot,\cdot) \mid  x\in N_1\}=\{D^2_y(\phi(\cdot),\phi(\cdot))\mid y \in N_2\}.
$
\end{enumerate}

\begin{theorem}[Distance difference functions determine the manifold \cite{LSaksala}]\label{thm dist diff}
Let $(N_1,h_1)$ and $(N_2,h_2)$ be compact and connected Riemannian manifolds of dimension $n\geq 2$ without boundary. 
Let $F_j\subset N_j$, $j=1,2$, be closures of non-empty open sets having a smooth manifold. Suppose the above conditions
(i) and (ii)  
are valid. 
 {Then the
manifolds $(N_1,h_1)$ and $(N_2,h_2)$ are isometric.}
\label{main theorem LS}
\end{theorem}

We will not consider here the proof of Theorem \ref{thm dist diff} which is based on Riemannian 
geometry. Instead, we discuss next how the distance difference functions  on a Riemannian manifold 
$(N,h)$ determine the
earliest light observation sets on the Lorentzian manifold $M=\R\times N$ with the product Lorentzian
metric $g=-dt^2+h$ and consider the determination of the Lorentzian metric using the knowledge of these sets. 
Let $(N,h)$  be a Riemannian manifold. We denote the unit speed geodesic of $(N,h)$ emanating from a point $(y,\eta)\in SN$  by $\alpha^h_{y,\eta}(t)=\exp_p(t\eta)$. 
Here, $SN=\{(y,\eta)\in TN;\ \|\eta\|_h=1\}$ is the unit sphere bundle of $N$. Using the Riemannian manifold $(N,h)$   we define
a product  manifold $M=\R\times N$. When $U\subset N$ is an open set having coordinates $X:U\to \R^n$,
$X(x)=(x^1(x),x^2(x),\dots,x^n(x))$ we define a Lorentzian metric $g=-dt^2+h(x)$, that is given in local coordinates $(t,x^1,\dots,x^n)\in \R\times U$ by
$$
 g=-dt^2+\sum_{j,k=1}^n h_{jk}(x^1,\dots,x^n)dx^jdx^k.
$$
Below, we denote the time coordinate also by $x^0=t$. Let $p=(t_0,y)\in \R\times N$ and $\xi=\frac \p{\p t}+\eta$ be a vector at the point $p$, where $(y,\eta)\in SN$.
Observe that $g(\xi,\xi)=-1+\|\eta\|_h^2=0$, that is, $\xi$ is a future pointing light-like vector.
 The curve $\gamma_{p,\xi}:\R\to \R\times N$,
 \ba
\gamma_{p,\xi}(s)=(t_0+s,\alpha^h_{y,\eta}(s))
 \ea
 is a geodesic of the Lorentzian manifold $(M,g)$.  Also, let ${{U}}\subset \R\times F$ be an open set where we consider observations. For the point $p\in M$  we define the light-observation set with the observation set 
 ${{U}}\subset M$ by setting
 \ba
 \mathcal P_{{U}}(p)=
 \{\gamma_{p,\xi}(s)\in M\mid  s\ge 0,\  \xi=\frac \p{\p t}+\eta,\ \|\eta\|_h=1\}\cap {{U}}.
 \ea
 Let $\kappa:SN\to \R\cup \{\infty\}$  be the cut-locus function on the manifold $(N,h)$, defined by
 \ba
 \kappa(y,\eta)=\sup\{ s\ge 0:\ \dist_h(\alpha_{y,\eta}(s),y)=s\}
 \ea
 Then, the earliest light-observation set of the point $p=(t_0,y)$ with the observation set 
 ${{U}}\subset M$   is
 \ba
 \mathcal E_{{U}}(p)=\{\gamma_{p,\xi}(s)\in M\mid 0\le s\le \kappa(y,\eta),\  \xi=\frac \p{\p t}+\eta,\ \|\eta\|_h=1\}\cap {{U}}.
 \ea
Next, we consider the earliest light-observation sets on general Lorentzian manifolds,  which may not be of the product type.

\subsection{Inverse problems on Lorentzian manifolds}
\subsubsection{Notations in a Lorentzian space-time}
\label{sec1}

Let $(M, g)$ be a $(1+n)$-dimensional time-oriented Lorentzian manifold of signature $(-,+,+,\dots,+)$, where $n\ge 2$. On $M$, we use local coordinates $(x^0,x^1,\dots,x^n)$, where $x^0=t$
it the time-coordinate. 

\modified{Let $L_pM=\{\xi\in T_pM\setminus\{0\}:\ g(\xi,\xi)=0\}$ {\itext be the set of light-like vectors 
in the tangent space $T_pM$}.

Let $\exp_p: T_pM \to M$
be the exponential map on $(M,g)$.
The geodesic starting at $p$ in the
direction $\xi \in T_pM\setminus\{0\}$ is the curve
 $\g_{p, \xi}(t)=\exp_p(\xi t)$,  $t \geq 0$.  

Let 
$\mu:[-1,1]\to M$ be a smooth future pointing time-like geodesic
and ${{U}}\subset M$ be an open neighborhood of $\mu([-1,1])$.

Let $q \in M$.
The set of future 
pointing \emph{light-like} vectors at $q$ is denoted by
\begin{equation*}
{L}_q^+ M = \{ \theta \in T_q M\setminus 0: \ g(\theta,\theta) = 0,\ \text{$\theta$ is future-pointing}\}.
\end{equation*}
A vector $\theta \in T_q M$ is \emph{time-like} if $g(\theta,\theta) < 0$ and \emph{space-like} if $g(\theta,\theta) > 0$.
\emph{Causal vectors} are the collection of time-like and light-like vectors, and a curve $\gamma$ is {time-like} (light-like, causal, future-pointing) if the tangent vectors $\dot \gamma$ are time-like (light-like, causal, future-pointing).

For $p, q \in M$, {\ftext the notation} $p \ll q$ means that $p, q$ can be joined by  a {\ftext future-pointing} time-like curve. 
The \emph{chronological future} and \emph{past} of $p \in M$ are
\begin{equation*}
I^+(p) = \{q \in M\ : \ p \ll q\}, 
\quad
I^-(p) = \{q \in M\ : \ q \ll p\}.
\end{equation*}
%
Similarly, for $p, q \in M$, {\ftext the notation} $p < q$ means that $p\not =q$ and $p, q$ can be joined by  a {\ftext future-pointing} causal (i.e., time-like or light-like) curve. 
We denote $p\le q$ if $p < q$ or $p=q$.
The \emph{causal future} and \emph{past} of $p \in M$ are
\begin{equation*}
J^+(p) = \{q \in M\ : \ p \le q\}, 
\quad
J^-(p) = \{q \in M\ : \ q \le p\}.
\end{equation*}
We will  denote throughout the paper the causal diamond-sets
\begin{equation}
\label{causal diamond}
I(p,q) = I^+(p) \cap I^-(q),\quad J(p,q) = J^+(p) \cap J^-(q).
\end{equation}
To emphasise the Lorentzian structure of $(M,g)$ we sometimes write 
 $I^\pm_{M,g}(p) = I^\pm(p)$ and  $J^\pm_{M,g}(p) = J^\pm(p)$.

A time-oriented Lorentzian manifold $(M,g)$ is \emph{globally hyperbolic} if
there are no closed causal paths in $M$, and for any $p,q \in M$ the set $J(p,q)$ is compact \cite{Bernal}.
{\ftext According to}  \cite{Bernal2}, a globally hyperbolic manifold is isometric to the product manifold $\mathbb{R} \times N$ with the Lorentzian metric given by
\begin{equation}\label{special metric}
g = -\beta(t,y) dt^2 + \kappa(t,y),
\end{equation} where $\beta : \mathbb{R} \times N \rightarrow \mathbb{R}_+$ and $\kappa$ is a Riemannian metric on $N$ depending on $t$.

\subsection{Passive measurements on a Lorentzian manifold}

%
%
%

\subsubsection{The set of earliest light observations}

{\revtext Let $s_-,s_+\in (-1,1)$, $s_-<s_+$,  and
$p^\pm=\mu(s_\pm)$. 
Next define the light observation sets and consider in particular the case when
 $W \subset I^-(p^+) \setminus J^-(p^-)$ is a relatively compact open set,
 see Fig.\ 2 (left).


\begin{definition} \label{def. O_U}
Let $(M,g)$ be a globally hyperbolic Lorentzian manifold.

(i) For $q\in M$, let
$$\L^+(q)=\exp_q(L^+_qM)\cup\{q\}=
\{\gamma_{q,\xi}(s)\in M: \ \xi\in L^+_qM,\ s\geq 0\}\subset M
$$ {\itext be  the future directed light-cone emanating from the point $q$.

The light observation set of $q$ in the observation set ${{U}}$} is 
$
\P_{{U}}(q)= {\mathcal L}^+(q) \cap {{U}} \in 2^{{U}}.
$

(ii) The earliest light observation set of $q\in M$  in ${{U}}$ is
\beq\label{eq BB a}
& &\hspace{5mm}\firstpoint_{{U}}(q)=\{x\in \P_{{U}}(q):\ \hbox{there exist no point $y\in\P_{{U}}(q)$ and}\\
\nonumber
& &\hspace{25mm}\hbox{
a future-pointing time-like path 
$\alpha:[0,1]\to {{U}}$}\\
\nonumber
& &\hspace{25mm}\hbox{ such that $\a(0)=y$ and $\a(1)=x$}\}\subset {{U}}
\eeq

(iii) 
Let $W\subset M$ be open.
 The collection of the earliest light observation  sets with source points at $W$ is 
\beq
\firstpoint_{{U}}(W)=\{\firstpoint_{{U}}(q): \ q\in W\}\subset 2^{{U}}.\label{collection}
\eeq
Note $\firstpoint_{{U}}(W)$
is defined as an unindexed set, that is, for an element $\firstpoint_{{U}}(q)\in \firstpoint_{{U}}(W)$ we do not
know what is the corresponding point $q$.
\end{definition}
Above, $2^{{U}}=\{{{U}}^\prime:\ {{U}}^\prime \subset {{U}}\}$ is the power set of ${{U}}$. 
{\gtext Below, when we say that the set ${{U}}$ is given as a differentiable manifold,
we mean that we are given the set ${{U}}$ and the local coordinate charts on it for which
the corresponding transition maps are $C^\infty$-smooth.}

%
%

The set $\mathcal{P}_{{U}}(q)$ can be viewed as a model of a measurement where light emitted by a point source at $q$ is recorded in ${{U}}$.
As gravitational wave packets propagate at the speed of light, 
$\mathcal{P}_{{U}}(q)$ could also correspond to an observation 
where a gravitational wave is generated at $q$
and detected in ${{U}}$.

\begin{figure}
\centering
{\includegraphics[height=3.5cm]{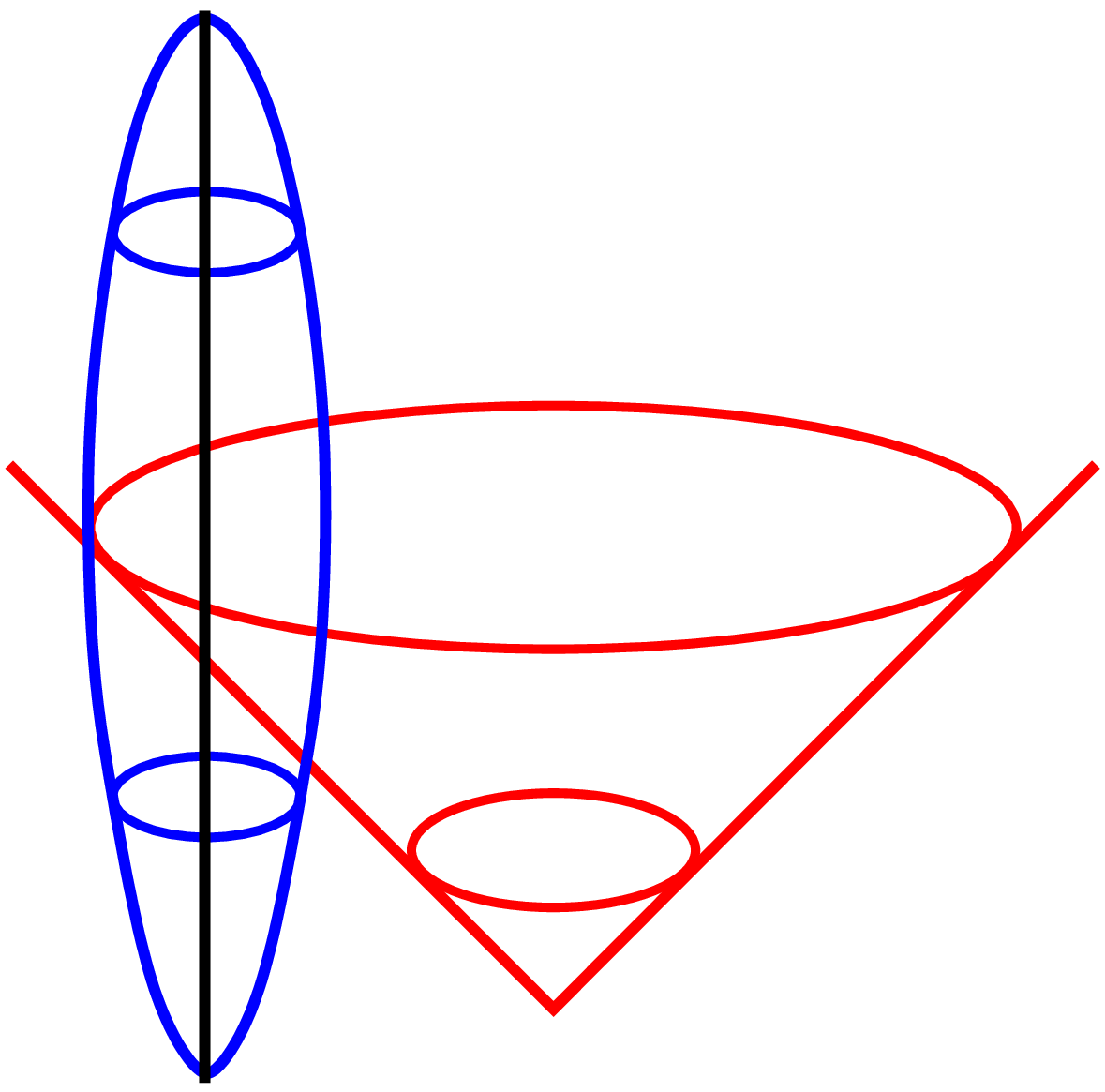}\quad\quad\quad
\includegraphics[height=3.5cm]{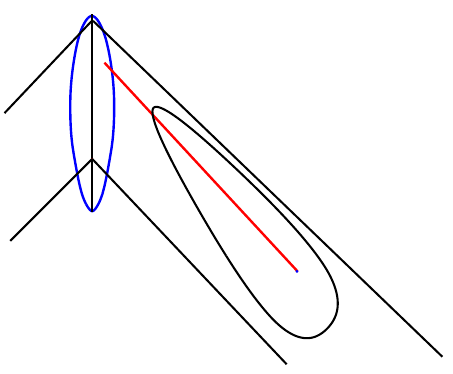}}
\caption{
{Left.:} 
When there are {no cut points},
the earliest light observation set $\mathcal{E}_{{U}}(q)$ is the intersection of the cone and the  open set  ${{U}}$. The cone is the union of 
 future-pointing light-like geodesics from $q$, and the ellipsoid depicts ${{U}}$.  Right: The unknown set 
 $W \subset I^-(\mu(1)) \setminus I^-(\mu(-1))$
 contains source points. Light cones emanating from these
 points (red segment) are observed in the set $U$.
 }
\end{figure}

The following theorem says, roughly speaking, that observations of a large number of point sources in a region $W$ determine the structure of the spacetime in $W$, up to a conformal factor.

\begin{theorem}[\cite{kurylev2018inverse}] \label{thm3}
Let $(M_j,g_j)$, {\ftext where} $j = 1,2$, be two open globally hyperbolic Lorentzian manifolds of dimension $1+n$, $n\geq 2$. 
Let $\mu_j : [-1,1] \to M_j$ be a  {\ftext future-pointing} time-like path, let ${{U}}_j \subset M_j$ be a neighbourhood of $\mu_j([-1,1])$,
and let  $W_j\subset M_j$ be such that
$W_j \subset I_{M_j, g_j}^-(\mu_j(1)) \setminus I_{M_j, g_j}^-(\mu_j(-1))$
are open and relatively compact sets.
Assume that there is a conformal diffeomorphism
$\phi:{{U}}_1 \rightarrow {{U}}_2$
such that $\phi(\mu_1(s)) = \mu_2(s)$, $s \in [0,1]$, and
\begin{equation*}
\{\phi(\mathcal{E}_{{{U}}_1}(q)):\ q \in W_1\} 
= \{\mathcal{E}_{{{U}}_2}(q) :\ q \in W_2\}.
\end{equation*}
Then there is a diffeomorphism
$\Psi : W_1 \rightarrow W_2$
and a strictly positive function $\alpha \in C^\infty(W_1)$
such that $\Psi^* g_2 = \alpha g_1$ and $\Psi \rvert_{{{U}}_1 \cap W_1} = \phi$.
\end{theorem}


Next, we consider an alternative formulation of the above inverse problem  where the data are given in terms of functions, not as a collection of sets.

Assume that  open set  ${{U}}\subset M$ where observations are made, is a union of time-like curves $\mu_a(-1,1)$, $a\in A$, where
$\mu_a:(-1,1)\to M$, $a\in  \overline A\subset \R^k$.
that depend continuously on $a$.
  
   Let  $p_1,p_2\in \mu_{a_0}$, $a_0\in A$. Also, let $W\subset J^-(p_2)\setminus J^-(p_1)$ be an  open,
relatively compact set.
The {observation time function} $F_q:\overline A\to \R$ for a point $q\in W$ is 
\ba
F_q(a)=\min\{s\in \R&:&\hbox{there exists a future-directed light-like }
\\ 
& &\hbox{geodesic from $q$ to $\mu_a(s)\}$}.
\ea 
That is, $F_q(a)$  is the first time when we observe on $\mu_a$ light that comes from $q$.

By Def.\  \ref{def. O_U}, the set
 $\E_U(q)$ determines the the function $F_q$ in $A$  in via the formula
\beq\label{e1 repeated }
\MTEXT{F_q(a)=s,\ \hbox{where $s\in[-1,1]$ is such that } \mu_a(s) \in \E_U(q),\  a\in {A}. \hspace{-.5cm}}
\eeq
As $F_q:\overline {A}\to \R$ is continuous, this determines $F_q(a)$ for all $a\in \overline{A}$.}
On the other hand,  $F_q$ determines $\E_U(q)$ via the formula
 \beq\label{eq: From F to e}
 \E_U(q) =\{ \mu_a(F_q(a)): \ a\in {A}\}.
 \eeq
Due to this, Theorem \ref{thm3} has the following equivalent formulation:

%
%
%
%
%
%
%

%
%
%
%
%
%

\begin{theorem}[\cite{kurylev2018inverse}]\label{main thm}
Let $(M,g)$ be a globally hyperbolic 
 Lorentzian manifold of dimension  $1+n\geq 3$.
Assume that  $\mu_a(-1,1)\subset M$, $a\in {A}\subset \R^m$ are  time-like
paths, $ {{U}}=\bigcup_{a\in A}\, \mu_a$ is open, and $p_1,p_2\in \mu_{a_0}$.
   Let  $W\subset J^-(p_2)\setminus J^-(p_1)$ be a relatively compact open
set and $W'\subset W$ be dense. 
Then  $({{U}},g|_{{U}})$ and the collection of the observation time functions,
 \ba
\big\{\ F_q:A\to \R \ \big|\ \  q\in W'\big \}\subset C({A}),
\ea
determine  the topological and differentiable type of the manifold
$W$ and the conformal class of the metric $g$ in  $W$.
\end{theorem}

We will outline the ideas of the proof of Theorem \ref{main thm} below in Section 
\ref{section proofs}.

\subsection{Active measurements: Inverse problems for non-linear wave equations}
{Let $(M,g)$ be a $(1+n)$-dimensional  globally hyperbolic Lorentzian manifold 
and {\ftext assume, without loss of generality, that $M = \R \times N$ with a metric of the form (\ref{special metric}).}
Let $t_0 > 0$ and consider the semilinear wave equation
\begin{eqnarray} \label{eq8} 
&\square_g u(x) + a(x) u(x)^m = f(x), \quad &\text{for $x\in (-\infty, t_0) \times N$},
\\
\label{init_cond_wave}
&u = 0, \quad &\text{in $(-\infty, 0) \times N$},
\end{eqnarray} 
where $m \ge 2$ is an integer.
Here $a \in C^\infty(M)$ is {\ftext a} nowhere vanishing function and $\square_g$ is the wave operator, or the Lorentzian Laplace operator,
\begin{equation*}
	\square_g u= \sum_{j,k=0}^n |\det(g)|^{-1/2} \frac{\partial}{\partial x^j}\left(|\det(g)|^{1/2} g^{jk} \frac{\partial}{\partial x^k} u\right).
\end{equation*}
Let $\mu: [-1,1] \to M$ be a time-like curve and ${{U}}$ be its open, relatively compact neighbourhood.
The solution of {\ftext(\ref{eq8}{\ftext )--(}\ref{init_cond_wave})} exists when the source $f$ is supported in ${{U}}$
and satisfies $\|f\|_{C^k(\overline {{U}})}<\e$, where $k\in \Z_+$  is sufficiently large and $\e>0$ is sufficiently small.
For such sources $f$ we define the measurement operator
\begin{equation}
L_{{U}} : f \mapsto u \rvert_{{U}}.
\end{equation} 

\begin{theorem}[\cite{kurylev2018inverse,FeizLO}] \label{thm4}
Let $(M,g)$ be a globally hyperbolic {\ftext $(n+1)$-dimensional} Lorentzian manifold, $n\ge 2$. Let $\mu$ be a time-like path containing $p^+$ and $p^-$. Let ${{U}} \subset M$ be a neighborhood of $\mu$ and {\ftext let} $a : M \rightarrow \mathbb{R}$ be a nowhere vanishing $C^\infty$-smooth function. Then $({{U}},g \rvert_{{U}})$ and the measurement operator $L_{{U}}$ determine the topology, differentiable structure and the conformal class of the metric $g$ in the double cone  $I_{M,g}(p^-, p^+)$.
\end{theorem}

We will outline the ideas of the proof of Theorem \ref{thm4} below in Section 
\ref{section proofs}.

%

\begin{figure}
\begin{center}

\arxivpreprint{
$  $\hspace{-10cm}
\includegraphics[height=3.5cm]{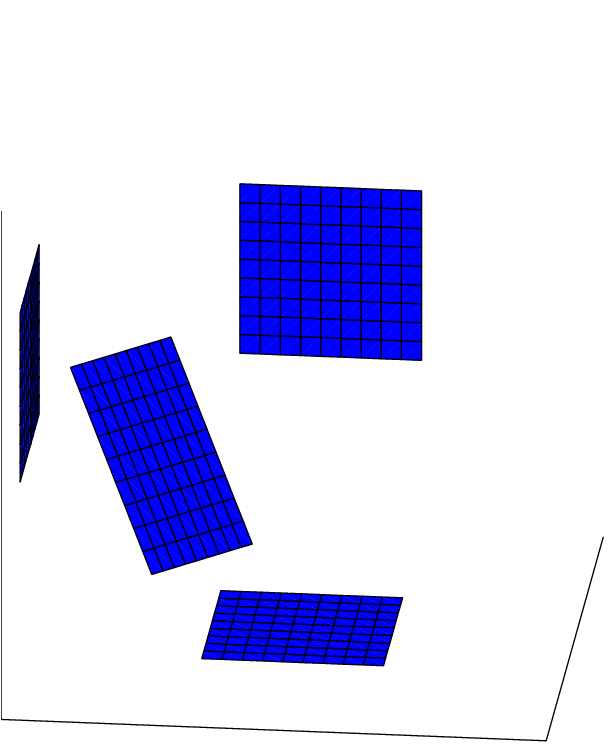}\quad 
\includegraphics[height=3.5cm]{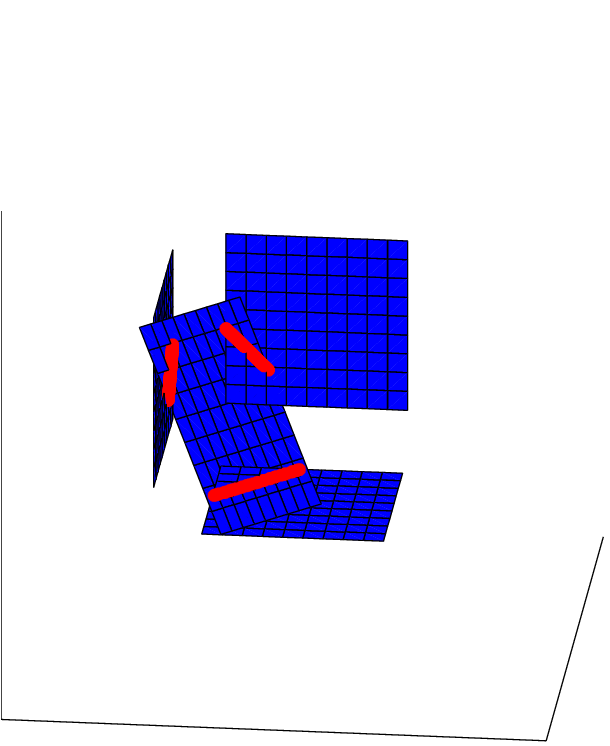}\quad
\includegraphics[height=3.5cm]{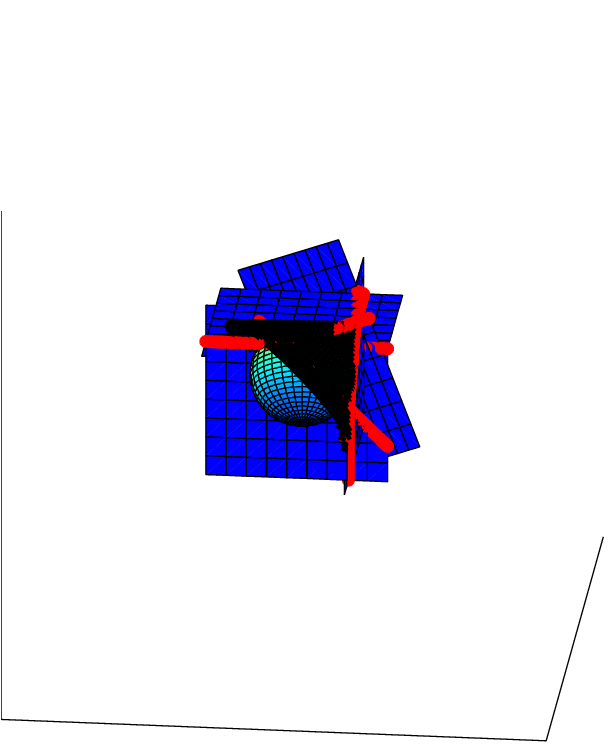}\quad
\includegraphics[height=3.5cm]{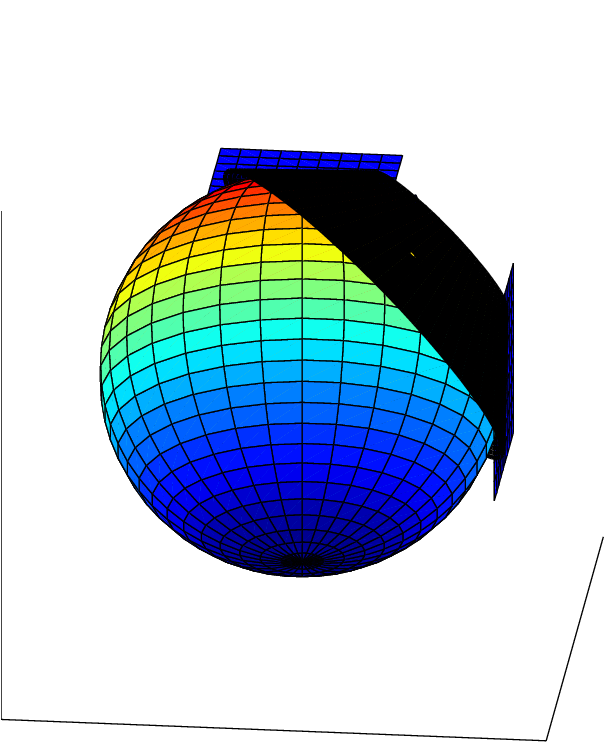}\hspace{-10cm}$  $
}
\end{center}
\caption{Four plane waves propagate in space. 
When the planes intersect, the non-linearity of the hyperbolic system  produces
new waves.
{\it Left:} Plane waves before  interacting.
 {\it Middle left:}
 The {\ftext two}-wave interactions (red line segments) appear  but do not cause
 singularities propagating to new directions. {\it Middle right and Right:}  All plane waves have intersected and new 
 waves have appeared. The  {\ftext three-}wave interactions cause 
conic waves ({\ftext the} black surface). Only one such wave is shown in the figure. The {\vtext interaction of four waves} causes
a microlocal point source  that sends a spherical wave in all future light-like 
directions. } 
\end{figure}

The proof of Theorem \ref{thm4}  {\ftext uses the results on} the inverse problem for passive measurements for point sources, described below,
and the 
 non-linear interaction of
 waves having conormal singularities.
There are many results on such non-linear interaction,  starting {\ftext with} the studies of 
 Bony \cite{Bony} and Melrose and Ritter \cite{MR1}.
%

}

\subsubsection{Stability results and numerical reconstruction}
As discussed below, any globally hyperbolic Lorentzian manifold $M$ is isometric to a product manifold $\R\times N$ equipped with the product metric \eqref{special metric}.
Let
$\Omega \subset N$ be a smooth submanifold of dimension $n$ with smooth boundary and  denote the lateral boundary of $[0,T]\times \Omega\subset M$ by 
$
\Sigma := [0,T]\times \p\Omega.
$
We denote the Sobolev space of functions having $s$ weak
derivatives on $\Sigma$ by $H^{s}(\Sigma)$ and
the closure of $C^\infty_0(\Sigma)$ in $H^{s}(\Sigma)$
is denoted by $H_0^{s}(\Sigma)$.

We consider the nonlinear wave equation
\begin{equation}\label{eq:Main equation}
\begin{cases}
\square_g u(t,x) + q(t,x)u(t,x)^m=0,&\text{in } [0,T]\times \Omega,\\
u=f,&\text{on } [0,T]\times \p \Omega,\\
u(0,x)=\p_t u(0,x) = 0,&\text{on } \Omega,
\end{cases}
\end{equation}
where we assume that the exponent $m$ is an integer greater or equal than $4$.

Let $\W$ be a compact set belonging to both the {chronological} future $I^+(\Sigma)$ and past of $I^-(\Sigma)$ the lateral boundary $\Sigma=[0,T]\times \p\Omega$:
\begin{equation}\label{eq:recovery_set}
\W\subset I^-(\Sigma) \cap I^+(\Sigma)\cap ([0,T] \times \Omega).
\end{equation}
\begin{theorem}[Stability estimate \cite{TyniLorentzian}]\label{thm:stability}
Suppose $(M,g)$, $M=\R\times N$, is a $(1+n)$-dimensional globally hyperbolic Lorentzian manifold.

Let $T>0$ and let $\Omega\subset N$ be an open set with smooth non-empty boundary. Let $m\geq 4$ be an integer, $s\in\N$ with $s+1>(n+1)/2$ and $r\in \R$ with $r\leq s$. Let $j=1,2$. Assume that $q_j\in C^{s+1}(\R\times\Omega )$ satisfy $\Vert q_j \Vert_{C^{s+1}} \leq c_0$, $j=1,2$, for some $c_0>0$. Let $\Lambda_j :H_0^{s+1}(\Sigma)\to H^r(\Sigma)$ be the corresponding Dirichlet-to-Neumann maps of the nonlinear wave equation~\eqref{eq:Main equation}.

Let $\eps_0>0$, $L>0$ and $\delta\in (0,L)$ be such that
 \begin{equation}
 \Vert \Lambda_1(f)-\Lambda_2 (f)\Vert_{H^r(\Sigma)} \leq \delta
\end{equation}
for all $f\in H_0^{s+1}(\Sigma)$ with $\Vert f\Vert_{H^{s+1}(\Sigma)}\leq \eps_0$. 
Then there exists a constant $C>0$, independent of $q_1,q_2$ and $\delta>0$, such that
\begin{equation}\label{eq:esimate_for_potential_difference}
  \norm{q_1-q_2}_{L^\infty( \W )}\leq C \delta^{\sigma(s,m)},
\end{equation}
where $\sigma(s,m) = \frac{8(m-1)}{2m(m-1)(8s-n+13)+2m-1}$.
\end{theorem}

Numerical reconstructions based on Theorem  \ref{thm:stability} in the $(1+1)$-dimensional case
are shown in Figure \ref{fig:EX4}.

\begin{figure}[tp]
\centering
\includegraphics[width = 0.8\textwidth]{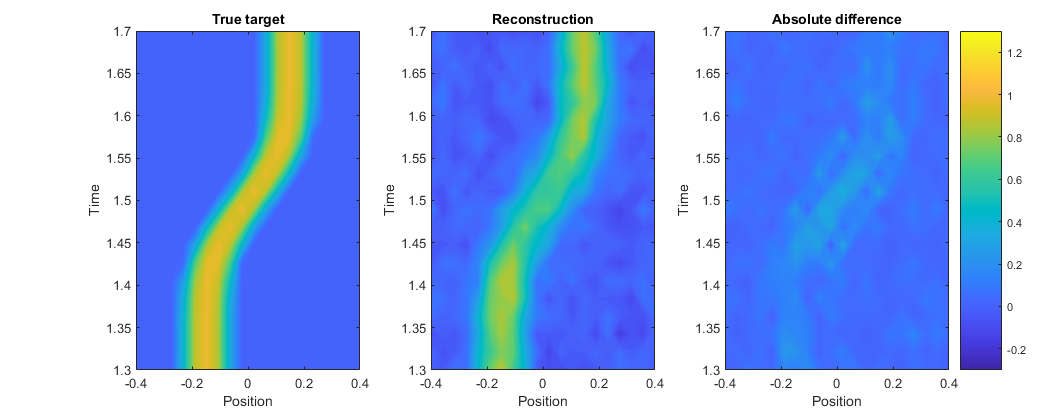}

\caption{Reconstruction of the coefficient $q(x,t)$ in
equation  $\p_t^2u(x,t)-\p_x^2u(x,t) + q(x,t) u(x,t)^{2} =0$, $ x\in [0,1],$ $t\in [0,T],$ from the data $f\to \bra \psi_0,\Lambda_q f\cet_{L^2({\p \Omega\times [0,T]})}$, where $\Lambda_q$ is the Dirichlet-to-Neumann map $\psi_0$ is a suitable chosen `device function'. Above, data has additive Gaussian noise and the signal-to-noise ratio is 12 dB.
Left figure: Ground truth of the potential $q(x,t)$. Center figure: Numerical reconstruction. Right figure: Error in reconstruction. For numerical results, see \cite{Tyninumerical} and for mathematical analysis, see \cite{TyniStability,TyniLorentzian}.
\label{fig:EX4}
}
\end{figure}

\section{Ideas for proofs and reconstruction methods}
\label{section proofs}

\subsubsection{Light-observation sets and the inverse problem for passive observations}

In this subsection, we outline the proof of Theorem \ref{main thm}. 
Below, let $\overline W=\hbox{cl}\,(W)$ be the closure of $W$ in $M$,
such that $\overline W\subset  J^-(p^+)\setminus I^-(p^-)$.
{The set $\overline {A}\times \overline W$ is compact and
the map $ \F:\overline W\to C(\overline {A})$, $\F(q)=F_q$, is continuous. 

Next we will show that the map 
$\F:\overline W \to \F(\overline W)$, $\F(q)=F_q$    is injective. 
To prove this, we assume the opposite: Assume that there are   $q_1,q_2\in W$ that satisfy $F_{q_1}=F_{q_2}$.
Then,
\beq\label{E sets are the same}
 \E_{{U}}(q_1)= \E_{{U}}(q_2).
 \eeq 
This implies that all light-like geodesics from $q_1$ to ${{U}}$ go through $q_2$ and hence, $q_1$  is in the past of $q_2$, or
vice versa. Without loss of generality, we can assume that  $q_1$  is in the past of $q_2$. Let $a_1\in {A}$
and $p_1=\mu_{a_1}(f^+_{q_1}(a_1))\in \mu_{a_1}([-1,1])$.
Let now $\gamma_{q_1,\xi_1}([0,\ell_1])$  be a light-like geodesics from $q_1$  to $p_1$, that is, $p_1=\gamma_{q_1,\xi_1}(\ell_1)\in \mu_{a_1}([-1,1])$.
Next, let  $a_2\in {A}$  be such that 
$\gamma_{q_1,\xi_1} \cap \mu_{a_2}([-1,1])=\emptyset$. Let
 $p_2=\mu_{a_2}(f^+_{q_1}(a_2))\in \mu_{a_2}([-1,1])$
and $\gamma_{q_1,\xi_2}([0,\ell_2])$  be a light-like geodesics from $q_1$  to $p_2$, that is, $p_2=\gamma_{q_1,\xi_2}(\ell_2)\in \mu_{a_2}([-1,1])$. Then, as $\gamma_{q_1,\xi_1} \cap\mu_{a_2}([-1,1])=\emptyset$, $\xi_2$  and $\xi_1$ are not parallel.
As all light-like geodesics from $q_1$ to ${{U}}$ go through $q_2$, we see that there are $t_1\in (0,\ell_1)$ and 
$t_2\in (0,\ell_2)$ such that 
$$
q_2=\gamma_{q_1,\xi_1}(t_1),\quad q_2=\gamma_{q_1,\xi_2}(t_2).
$$
As the vectors $\xi_2$  and $\xi_1$ are not parallel, the geodesics $\gamma_{q_1,\xi_1}$  and $\gamma_{q_2,\xi_2}$  are 
not the same paths, and we see that the path $\beta =\gamma_{q_1,\xi_1}([0,t_1])\cup \gamma_{q_1,\xi_2}([t_2,\ell_2])$
is a union of two causal paths from the point $q_1$ to $q_2$  and then from $q_2$  to $p_2$.
Since the paths the geodesics $\gamma_{q_1,\xi_1}$  and $\gamma_{q_2,\xi_2}$  are not the same paths with
different parametrizations, the path $\beta$   is not a geodesic. As the manifold $(M,g)$ is globally
hyperbolic, this implies 
that there is a time-like path from $q_1$ to $p_2$ and hence it is not possible 
 $p_2=\mu_{a_2}(f^+_{q_1}(a_2)\in \mu_{a_2}([-1,1])$, as we assumed above.
 This contradiction implies   $F_{q_1}\not =F_{q_2}$.
Thus, the map 
$\F:\overline W \to \F(\overline W)$, $\F(q)=F_q$    is injective.


The map $\mathcal F:\overline  W\mapsto C(\overline {{A}})$ is continuous and injective.
As $\overline W$ is compact, these imply that the map $\mathcal F:\overline  W\mapsto C(\overline {{A}})$ is a homeomorphism.}
Thus, $\mathcal F(\overline  W)\subset C(\overline {{A}})$ and $\overline W$ are homeomorphic.
Furthermore, as $W$ is an open subset of the manifold $M$, we can see that $\mathcal F( W)\subset C(\overline {{A}})$ and $W$ are homeomorphic,
so we can consider $ {\F}(  W)$  as a topological manifold which is a
homeomorphic copy of the original manifold $W$. Thus, we have reconstructed $W$ as 
a topological manifold.


 \begin{figure}[tp]


\centering{\includegraphics[height=3cm]{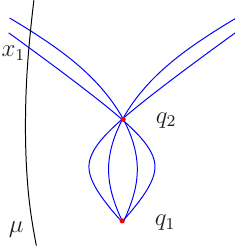}\quad \includegraphics[height=3cm]{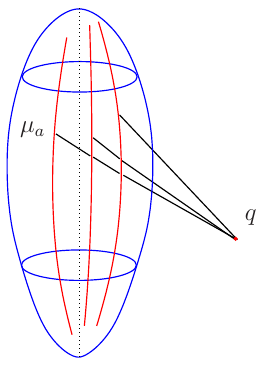}}



 \caption{Left: Picture on proving that $q_1=q_2$ by contradiction. Right: The function $X_{\vec a}:q\to (F_q({a_j})_{j=1}^{n+1}=(f^+_{a_j}(q))_{j=1}^{n+1}\in \R^{n+1}$ defines coordinates near a point $q$ on $M$.}
 \end{figure}

  Next,  we consider the reconstruction of $W$ as 
a differentiable manifold. Below, we will only consider how the
coordinate functions on $ {\F}(W)$ are constructed and omit the proof
that this constructions gives a differentiable structure on  $ {\F}(W)$
which makes it diffeomorphic to $W$. For details on this, we cite to \cite{kurylev2018inverse}.

For $a\in {A}$, let us define the evaluation function $E_a[F]=F(a)$. 
 For $q_0\in W$, let us choose $n+1$ parameters $a_j\in {A}$, $j=1,2,\dots,n+1$.
 When  the parameters $(a_1,a_2,\dots,a_{n+1})$ are chosen in a generic way
 that is, this $({n+1})$-tuple is in an open and dense subset of ${A}^{n+1}$, the function 
 $X_{\vec a}(q)=(F_{q}(a_j))_{j=1}^{n+1}=(f^+_{a_j}(q))_{j=1}^{n+1}\in \R^{n+1}$ defines smooth coordinates on $M$ near $q_0$ (Here, we omit this part of the proof, see \cite{kurylev2018inverse} on details).
 This implies that for any ${\F}(q_0)\in {\F}(W)$ and generic parameters $a_j\in {A}$, $j=1,2,\dots,{n+1}$
 there is a neighborhood
 $\mathcal V\subset   {\F}(W)$ such the function ${\bf X}_{\vec a}:F\to (E_{a_j}[F])_{j=1}^{n+1}\in \R^{n+1}$ defines coordinates in $\mathcal W$
such that the map $\F$ is a local diffeomorphism near the point $q_0$ from the manifold $W$ to 
$ {\F}(W)$.
 
Next, we will consider the construction of the Lorentzian metric on ${\F}(W)$ which
makes it conformal to $(W,g|_W)$. To this end, let us consider a point $f_0= {\F}(q_0)\in  {\F}(W)$
and its open neighborhood $\mathcal W\subset {\F}(W)$ where 
 $\vec a=(a_j)_{j=1}^{n+1}$ define coordinates ${\bf X}_{\vec a}:\mathcal W\to\R^{n+1}$ near $f_0$.

When $B\subset {A}$ is a sufficiently small neighborhood of the point $a_1$, the function $a\to  f^+_a({q_0})
$ is smooth  for $a\in B$. 
Let $d f^+_a$ denote the co-derivative of the map $f^+_a:W\to \R$.
We see that the co-vectors 
$
d f^+_a|_{q_0}\in T^*_{q_0}M,\quad a\in B
$
are light-like co-vectors at $q_0$. Let $G={\F}_*g|_W$ be the metric $g|_W$ pushed 
forward in the map ${\F}$ to the manifold ${\F}(W)$ and $E_a:{\F}(W)\to \R$
be the evaluation maps $E_a:f\to f(a).$  
Then, we see that the co-derivatives 
$$
dE_a|_{{\F}(q_0)}\in T^*_{{\F}(q_0)}{\F}(W),\quad a\in B
$$
are light-like co-vectors with respect to the metric $G$.

%
As $a\in B$  is above arbitrary, we see that the set  $\{s\,dE_a|_{{\F}(q_0)}\in T^*_{{\F}(q_0)}{\F}(W)\mid a\in B,\ s>0
\}$ contains an open subset of the light cone in the space $( T^*_{{\F}(q_0)}{\F}(W),G)$.
As the light cone  is a quadratic variety,
this set determines the whole light cone in $( T^*_{{\F}(q_0)}{\F}(W),G)$.
These sets determine the conformal class of the metric $G$ on ${\F}(W)$ at the point 
${\F}(q_0)$.

\subsubsection{Nonlinear interaction of waves and inverse problems for active measurements}

In this subsection, we outline the proof of Theorem \ref{thm4}.
To do this, we consider the  non-linear equation $\square_{}u+u^4=f$ in $\R^{1+3}$, that is, we  introduce the techniques by
considering the Minkowski space $\R^{1+3}$ and the wave operator
$\square=\p_{x^0}^2-\sum_{j=1}^3\p_{x^j}^2$.

Let us consider  the solutions $u_{\vec\e}(x)$  of  the non-linear wave equation
$$\square_{}u_{\vec\e}+(u_{\vec\e})^4=0,\quad
$$ 
that depend on parameters $\vec\e=(\e_1,\e_2,\e_3,\e_4)\in \R^4$.
Let us assume that $u_{\vec \e}|_{\vec \e=0}=0$. Then, the linearized waves
\ba
u_j(x)=\p_{\e_j}u_{\vec\e}\big|_{\vec\e=0},\quad j=1,2,3,4
\ea 
satisfy $\square_{}u_j=0$. Moreover, the 4th order derivative
\beq\label{4th order derivative}
w=\p_{\e_1}\p_{\e_2}\p_{\e_3}\p_{\e_4}u_{\vec \e}(x)\big|_{\vec\e=0}
\eeq
satisfies $$\square_{}w=S,\quad S=-24u_1u_2u_3u_4.$$
We say that $S$ is an artificial source produced by the 4-interaction of the waves.

 
We use coordinates $x=(t,y_1,y_2,y_3)\in \R^{1+3}$ and consider,
as an example, the
linearized waves $u_j(t,y_1,y_2,y_3)$ are 
\ba
& &u_1(t,y)=\delta(t-y_1),\quad
u_2(t,y)=\delta(t-y_2),\\
& &u_3(t,y)=\delta(t-y_3),\quad
u_4(t,y)=\delta(t-(y_1+y_2)/{\sqrt 2})
\ea
Then,
\ba
& &u_1u_2u_3u_4=\delta_{p}(t,y),\quad \hbox{where }p=(0,0,0,0) \in \R^{1+3}
\ea
and the 4-interaction of waves, $w(x)$, is the solution  of the wave equation
 \beq\label{eq S source}
 \square w=S,\quad S=-24u_1u_2u_3u_4=-24\delta_{p}.
 \eeq
In physical terms, $S$ corresponds to an (artificial) point source in the space time.

%

Let us next consider a general globally hyperbolic, $(1+3)$ dimensional Lorentzian manifold when 
$m=4$. 
The proof of Theorem \ref{thm4} is based on the analysis of {\ftext the} interaction of four waves.
Analogously to (\ref{4th order derivative}), the derivative $w=\partial_{\epsilon_1} \partial_{\epsilon_2} \partial_{\epsilon_3} \partial_{\epsilon_4} u_{\vec{\epsilon}} \big|_{\vec{\epsilon} = 0}$ of the solution $u_{\vec{\epsilon}}$
of (\ref{eq8}{\ftext )--(}\ref{init_cond_wave}) with the source
\begin{equation*}
f_{\vec{\epsilon}}(x) = \sum_{j=1}^4 \epsilon_j f_j(x), \quad \vec{\epsilon} = (\epsilon_1,\epsilon_2,\epsilon_3,\epsilon_4),
\end{equation*}
is  given by formula \eqref{eq S source},
where 
$u_j(x)$
are solutions to the wave equation $\square_g u = f_j$.
%

When $f_j$ is a $C^\infty$-smooth function outside the point $p_j\in M$, the linearized wave $u_j=\square_g^{-1}f_j$  produced
by the source $f_j$ is microlocally similar to Green's functions (i.e. the fundamental solution of 
the wave equation) corresponding to a point source at the point $p_j$. 
Let us denote the light cones in $M$ emanating from the point $p_j$ to future by $K_j$,
and, for simplicity, assume that the light-like geodesics of $M$ have no cut or conjugate points.
Then, in $M\setminus \{p_j\}$, the linearized solutions $u_j$ are  conormal distributions that are non-smooth on surfaces $K_j$ (see \cite{GU1,MU1}).
Then, two waves $u_1$ and $u_2$ are singular on hyper-surfaces
$K_1$ and $K_2$, respectively, and these singularities interact on $K_1\cap K_2$. 
This interaction creates only singularities that propagate along $K_1$ and $K_2$ on which the original singularities
 $u_1$ and $u_2$ propagate on.
The interaction of the three waves $u_1$, $u_2$, and $u_3$ happens on the intersection $K_{123}=K_1\cap K_2\cap K_3$ 
{\ftext which} is a curve in the space-time. As $N^*K_{123}$
contains light-like directions that are not in union, $N^*K_1\cup N^*K_2\cup N^*K_3$, this interaction produces interesting singularities that start to propagate. These singularities correspond to the black  conic wave in Fig.\ 3 (We note that these singularities produced by 3-interactions can also be used to determine the light-observation sets, see \cite{FeizLO}). 
Finally, singularities of all four waves $u_j$, $j=1,2,3,4$ interact at the point $\{q\}=\bigcap_{j=1}^4K_j$. 
The singularities from the point $q$ propagate along the light cone emanating from this point and with suitably chosen sources $f_j$  the wave  
$w$  is singular at any given point of the the light cone $\L(q)$.  Thus $S_{1234}$  can be considered as a microlocal point source that sends  singularities in all directions
as a point source located at the point $q$, and these singularities are observed in the set ${{U}}$.
 In this way,  the non-linear interaction
of waves determine the intersection of the light cone $\L(q)$ and the observation domain ${{U}}$.  Hence, we can determine the earliest light observation sets 
$\mathcal E_{{U}}(q)$ for an arbitrary point $q\in I(p^-,p^+)$. 
By letting $q$ vary through the set $I(p^-,p^+)$, we can apply Theorem \ref{thm3} to recover the topology, differentiable structure, and the conformal class of $g$ in $I(p^-,p^+)$.


\begin{funding}
The author was partially supported by the AdG project 101097198 of the European Research Council, Centre of Excellence of Research Council of Finland and the FAME flagship of the Research Council of Finland (grant 359186).
Views and opinions expressed are those of the authors only and do not necessarily reflect those of the funding agencies or the  EU.
\end{funding}

\def\cprime{$'$}

\end{document}